%
%

\magnification=1200


\font\AAA=cmr14 at 12pt
\font\BBB=cmr14 at 8pt

\overfullrule=0in

\def\Arr#1{\buildrel \hbox{#1}\over \longrightarrow}

\def\M{{\bf M}}
\def\oa#1{\overrightarrow #1}

\def\deg{{\rm deg}}

\def\log{{\rm log}}

\def\arr{\longrightarrow}
\def\supp{{\rm supp}\,}
\def\Link{{\rm Link}}
\def\Wind{{\rm Wind}}
\def\Div{{\rm Div}}


\def\Theorem#1{\medskip\noindent {\AAA T\BBB HEOREM \rm #1.}}
\def\Prop#1{\medskip\noindent {\AAA P\BBB ROPOSITION \rm  #1.}}
\def\Cor#1{\medskip\noindent {\AAA C\BBB OROLLARY \rm #1.}}
\def\Lemma#1{\medskip\noindent {\AAA L\BBB EMMA \rm  #1.}}
\def\Remark#1{\medskip\noindent {\AAA R\BBB EMARK \rm  #1.}}
\def\Note#1{\medskip\noindent {\AAA N\BBB OTE \rm  #1.}}
\def\Def#1{\medskip\noindent {\AAA D\BBB EFINITION \rm  #1.}}

\def\Conj#1{\medskip\noindent { \AAA C\BBB ONJECTURE \rm    #1.}}

\def\pf{\medskip\noindent {\bf Proof.}\ }
\def\qed{\hfill  $\vrule width5pt height5pt depth0pt$}
\def\equdef{\buildrel {\rm def} \over  =}

   \def\cc{{\cal C}}     
   \def\cp{{\cal P}}
   \def\co{{\cal O}}
\def\ce{{\cal E}}   \def\ck{{\cal K}}
\def\ch{{\cal H}}
\def\cs{{\cal S}}
\def\cd{{\cal D}}
\def\cl{{\cal L}}
\def\cp{{\cal P}}

\def\and{\qquad {\rm and} \qquad}
\def\arr{\longrightarrow}

\def\bbr{{\bf R}}
\def\bbc{{\bf C}}

\def\bbz{{\bf Z}}
\def\bbp{{\bf P}}

\def\a{\alpha}

\def\d{\delta}
\def\e{\epsilon}

\def\g{\gamma}

\def\l{\lambda}
\def\o{\omega}

\def\s{\sigma}
\def\x{\xi}

\def\D{\Delta}
\def\L{\Lambda}
\def\G{\Gamma}

\def\omp{quasi-plurisubharmonic }
\def\PH#1{\widehat {#1}}
\def\psh{{PSH}_{\o}}

\def\PL{3}
\def\PSH{4}
\def\PHC{5}
\def\PAW{6}

\def\PM{7}
\def\MH{8}

\def\g{\Gamma}

\  \vskip .2in
\centerline{\bf   PROJECTIVE LINKING  }
\smallskip

\centerline{\bf
 AND BOUNDARIES OF POSITIVE HOLOMORPHIC CHAINS }
\smallskip

\centerline{\bf IN PROJECTIVE MANIFOLDS, PART I  }
\vskip .2in
\centerline{\bf F. Reese Harvey and H. Blaine Lawson, Jr.$^*$}
\vglue
.5cm\smallbreak\footnote{}{ $ {} \sp{ *}{\rm Partially}$  supported by
the N.S.F. }

\centerline{\sl Dedicated to Nigel Hitchin }
\centerline{\sl in celebration of his 60th birthday}

 \vglue 1cm

\centerline{\bf Abstract} \medskip
  \font\abstractfont=cmr10 at 10 pt

{{\parindent=.5in\narrower\abstractfont \noindent 
We introduce the  notion of the {\sl projective linking number }
Link$_\bbp(\G,Z)$ of a compact oriented  real  submanifold $\G$ of dimension  $2p-1$ in complex projective $n$-space $\bbp^n$ with an algebraic subvariety 
$Z\subset \bbp^n-\G$ of codimension $p$. This notion is related to {\sl projective winding numbers}
and quasi-plurisubharmonic functions, and it  generalizes directly from  $\bbp^n$  to any  projective
manifold.  Part 1 of this paper establishes the following result for  the case $p=1$.
Let $\G$ be an oriented, stable, real analytic curve  in $\bbp^n$. Then \vskip .1in

\noindent{\sl
$\G$ is the boundary of a positive holomorphic
$1$-chain of mass $\leq \Lambda$ in $\bbp^n$ if and only if  
$\widetilde{\Link}_\bbp(\G,Z) \geq -\Lambda$ for all algebraic hypersurfaces 
$Z\subset \bbp^n-\G$.}

\vskip .1in\noindent
where $\widetilde{\Link}_\bbp(\G,Z) = {\Link}_\bbp(\G,Z)/\deg(Z)$.
An analogous theorem is implied in any  projective manifold.
Part 2 of this paper studies   similar  results for $p>1$.

 }}

 \vskip .35in

\centerline{\bf Table of contents.}
\vskip .25in

\moveright .6in\vbox{
1. Introduction.
 
2. Projective Hulls.

3. Projective Linking and Projective Winding Numbers.

4. Quasi-plurisubharmonic Functions.

5. Boundaries of Positive Holomorphic Chains

6. The Projective Alexander-Wermer Theorem for Curves.

7.  Theorems for General Projective  Manifolds.

8. Relative Holomorphic Cycles.
}

\vfill\eject

\centerline{\bf 1. Introduction}
\medskip

  In 1998 Herb Alexander and John Wermer published the following result [AW$_2$].
  
   \Theorem { (Alexander-Wermer)} {\sl  Let $\G$ be a compact oriented smooth submanifold of dimension 
  $2p-1$ in $\bbc^n$.  Then $\G$ bounds a positive holomorphic $p$-chain in $\bbc^n$ if 
  and only if  the linking number 
  $$
  {\rm Link}(\G, Z)\ \geq\ 0
  $$
  for all canonically oriented  algebraic subvarieties $Z$  of  codimension $p$ in $\bbc^n-\G$.  }

  \medskip

   The {\bf linking number} is an integer-valued  topological invariant defined by the intersection 
   ${\rm Link}(\G,Z) \equiv N\bullet Z$ with any 2$p$-chain $N$  having   $\partial N=\G$
   in  $\bbc^n$. (See \S 3.)  A {\bf positive holomorphic $p$-chain} is  a finite sum of canonically oriented complex subvarieties of dimension $p$ and finite volume in $\bbc^n -\G$.  (See Definition 5.1). Here the notion  of boundary is  taken in the sense of currents, i.e., Stokes' Theorem is satisfied. However,
   for smooth $\G$ there is boundary regularity almost everywhere, and if $\G$ is real analytic,
   one has complete boundary regularity. (See HL$_1$].)

  The main point of this paper  is to formulate and prove an analogue of the Alexander-Wermer Theorem               
  for oriented (not necessarily connected) curves  in a projective  manifold. In the sequel   we  shall study the corresponding result for  submanifolds of any odd dimension.

  Before stating the main result we remark that a key ingredient in the proof of the Alexander-Wermer Theorem 
  is the following  classical theorem and its generalizations [W$_1$].

 \Theorem{(Wermer) } {\sl Let   $\g\subset \bbc^n$ be a compact real analytic curve and denote by
  $$
  {\PH \g}_{\rm poly}   \ =\  \{z\in \bbc^n\ :\ |p(z)|\leq \sup_{\g} |p| \ \ {\sl for \ all \ polynomials}\ \  p\}.
  $$
  its polynomial hull. Then ${\PH \g}_{\rm poly}  - \g$ is a one-dimensional complex analytic subvariety of
  $\bbc^n-\g$.
  }
  \medskip
  
   For compact subsets $K$ of complex projective $n$-space $\bbp^n$, the authors recently introduced the notion of the {\bf projective hull} 
  $$
  {\PH K}_{\rm proj} \ =\ \{ z\in \bbp^n : \exists C  {\ \rm s.t.\ } \|\sigma_z\|\leq C^d \sup_K \|\sigma\| 
   \ \ \forall 
  \s\in H^0(\bbp^n, \co(d)),  d>0\}
  $$
  and defined  $K$ to be  {\bf stable} if  the constant $C$   can be chosen independently  of the point  $z\in  {\PH K}_{\rm proj} $.
   A number of basic properties  of  ${\PH K}_{\rm proj} $ were established  in [HL$_5$]),
  and  the following analogue of    Wermer's Theorem was proved in [HLW].

    \Theorem{(Harvey-Lawson-Wermer)} {\sl Let   $\g\subset \bbp^n$ be a stable   real analytic curve. Then ${\PH \g}_{\rm proj}  - \g$ is a one-dimensional complex analytic subvariety of
  $\bbp^n-\g$.
  }
  \medskip

  It is interesting to note that while Wermer's Theorem holds for curves with only weak differentiability 
  properties (see [AW$_1$] or [DL] for an account), its projective analogue fails even for $C^\infty$-curves.
  On the other hand there is much evidence for the following.

  \Conj{A} {\sl Every real analytic curve in $\bbp^n$ is stable.} \medskip

This brings us to the notion of projective linking numbers.  Suppose that $\g\subset \bbp^n$ 
  is a compact oriented  smooth curve, and let $Z\subset \bbp^n-\g$ be an algebraic subvariety of codimension $1$.  The {\bf projective linking number} of $\g$ with $Z$ is defined to be
 $$
  \Link_{\bbp}(\g, V) \  \equiv \ N\bullet Z - \deg (Z) \int_N\omega 
  $$
where $\omega$ is the standard  K\"ahler form on $\bbp^n$ and $N$ is any  integral 2-chain
with $\partial N=\g$ in $\bbp^n$. 
Here $Z$ is given the canonical orientation,
and $\bullet:H_2(\bbp^n,\g)\times H_{2n-2}(\bbp^n-\g)\to\bbz$ is the topologically defined 
intersection pairing.
This definition is independent of the choice of $N$. (See \S 3.)
  The associated {\bf reduced linking number} is defined to be
  $$
  \widetilde{\Link}_{\bbp}(\g,Z) \ \equiv \ {1\over \deg(Z)}\Link_{\bbp}(\g,Z)
  $$
  The basic result proved here is the following.
  
  \Theorem{\PAW.1} {\sl  Let $\g$ be a  oriented stable real analytic curve in $\bbp^n$ with a positive integer multiplicity on each component. Then the following are equivalent:
  \medskip
  
  (i)\ \ \ $\g$ is the boundary of a positive holomorphic $1$-chain of mass $\leq  \Lambda $ in  $\bbp^n$.
  \medskip
  
(ii)\ \  $ \widetilde{\Link}_{\bbp}(\g,Z)\ \geq \   -\Lambda$ for all algebraic hypersurfaces $Z$         in $\bbp^n-M$.
}
\medskip  

If $\g$ bounds any positive holomorphic 1-chain, then there is a unique such chain $T_0$ of least mass.
(All others are obtained by adding algebraic 1-cycles to  $T_0$.) Note that   $\Lambda_0 \equiv M(T_0)$ is the smallest positive number such that (ii) holds.

\Cor{\PAW.8} {\sl  Let $\g$ be as in Theorem \PAW.1 and suppose  $T$ is a
 positive holomorphic 1-chain with $dT=\g$.  Then $T$ is the unique holomorphic chain 
 of least mass with  $dT=\g$ if and only if  
 $$
 \inf_Z \left\{ {T\bullet Z \over \deg Z} \right\} \ =\ 0
 $$
where the infimum is taken over all algebraic hypersurfaces in the complement $\bbp^n-\g$.}
\medskip

Condition (ii) in Theorem \PAW.1 has several 
  equivalent  formulations.  The first is in terms of  {projective winding numbers}.  Given a holomorphic section  $\sigma \in H^0(\bbp^n, \co{}(d))$, the {\bf projective winding number}  of $\sigma$ on                  $\Gamma$ is defined as the integral
  $$
\Wind_{\bbp}(\Gamma, \s)\ \equiv\   \int_{\Gamma} d^C \log\|\s\|,
  $$
  and we set 
  $$
 \widetilde{ \Wind}_{\bbp}(\Gamma, \s)\ \equiv\  {1\over d}\Wind_{\bbp}(\Gamma, \s)
  $$
  
  Another formulation involves the cone  $\psh(\bbp^n)$ of {\bf \omp functions}. These are the upper semi-continuous functions $f:\bbp^n \to [-\infty, \infty)$ for which $dd^Cf+\o$ is a positive (1,1)-current on $\bbp^n$.
  
  \Prop{\PHC.2}  {\sl  Let $\Gamma$ be an oriented smooth curve in $\bbp^n$ with a positive integer multiplicity on each component. Then for any $\Lambda >0$ the following are equivalent:
 \medskip
 
 (ii)\ \ \  $ \widetilde{\Link}_{\bbp}(\G,Z)\ \geq \  -\Lambda$ for all algebraic hypersurfaces  
 $Z\subset\bbp^n-\G$.
 \medskip
 
 (iii) \ \ $ \widetilde{ \Wind}_{\bbp}(\Gamma, \s)\ \geq \ -\Lambda$ for all holomorphic sections $\s$ of 
 $\co{}(d)$,  and all $d>0$.
  \medskip
  
  (iv) \ \ $\int_\G d^C u\ \geq\ -\Lambda$  for all $u\in\psh(\bbp^n)$.
   }
  \medskip
  
  Any smooth curve $\Gamma$ in $\bbp^n$ lies in some affine chart $\bbc^n\subset\bbp^n$, and it is 
  natural to ask for a reformulation of condition (ii) in terms of the conventional linking numbers of
  $\G$ with algebraic hypersurfaces in that chart. This is done explicitly in Theorem \PAW.6.

The results above  extend to  any projective manifold
$X$.  Given a very ample line bundle $\lambda$ on $X$, there are intrinsically defined
{\bf $\lambda$-linking numbers} $\Link_{\lambda}(\Gamma, Z)$,  
{\bf $\lambda$-winding numbers} ${\Wind}_{\lambda}(\Gamma, \s)$ for $s\in H^0(X, \lambda^d)$,
and $\o=c_1(\lambda)$-{\bf quasi-plurisubharmonic functions} $\psh(X)$. With these notions, 
Proposition \PHC.2 and 
Theorem \PAW.1 carry over to $X$.
This is done in section \PM.  

  Theorem \PAW.1 leads to the following interesting result.
  
  \Theorem{8.1} {\sl Let  $\gamma\subset \bbp^n$  be  a finite disjoint union of
  real analytic curves and assume $\gamma$ is stable. Then a class
   $\tau\in H_2(\bbp^n,\gamma; \bbz)$ is represented by a positive holomorphic chain
   with boundary  on $\gamma$ 
   if and only if 
   $$
   \tau\bullet u\ \geq\ 0
   $$
  for all 
   $u\in H_{2n-2}( \bbp^n-\gamma;\bbz)$  represented by  positive algebraic hypersurfaces in 
   $\bbp^n-\gamma$.}
  \medskip
   
   This result expands to a   duality between the cones of relative and absolute classes
   which are representable by positive holomorphic chains [HL$_8$].
   
   The arguments given in \S 6 show that there is even more evidence for the following.
   
   \Conj{B} {\sl Every oriented real analytic curve in $\bbp^n$
   which satisfies the equivalent conditions of Proposition \PHC.2 is stable.} \medskip
    
  In Part II of this paper we prove that if Conjecture $B$ holds for curves in $\bbp^2$, 
  then all  the  results above continue to hold for  real analytic  $\G$ of any odd dimension $2p-1$.
  (No stability hypothesis is needed.)
 
  \Remark{}  There are several quite different characterizations of the boundaries
  of  general (i.e., not necessarily positive ) holomorphic chains  in projective and 
  certain quasi-projective manifolds.  See, for example, [Do], [DH$_{1,2}$],  and [HL$_{2,4}$].

\Note{}  To keep formulas simple throughout the paper  we adopt the
convention  that 
$$
d^C\ =\ {i\over 2\pi}(\overline\partial - \partial)
$$

\vfill\eject


\centerline{\bf 2. Projective Hulls}
\medskip
  
In this section we recall  the definition and basic properties of the projective
hull introduced in [HL$_5$]. This material is not really necessary for reading
the rest of the paper.

Let $\co(1)\to\bbp^n$ denote the holomorphic
line bundle of Chern class 1  endowed with the standard unitary-invariant
metric, and let $\co(d)$ be its $d$th tensor power with the induced tensor
product metric. 

\Def{2.1} Let $K\subset \bbp^n$ be a compact subset of complex projective
$n$-space. A point $x\in \bbp^n$ belongs to the {\sl projective hull of }
$K$ if there exists a constant $C=C(x)$ such that
$$
\|\s(x)\|\ \leq\ C^d \sup_K \|\s\| 
\eqno(2.1)
$$
for all  global sections $\s \in H^0(\bbp^n, \co(d))$ and all $d>0$.
This set of points is denoted $\PH K$.
\medskip

The set $\PH K$ is independent of the choice of metric on  $\co(1)$.

The projective hull possesses  interesting properties.  It and its
generalizations function in projective and K\'ahler manifolds much as the 
polynomial hull and its generalizations function in affine and Stein
manifolds. The following were established in [HL$_5$]:\medskip

\item {(1)} If $Y\subset \bbp^n$ is an algebraic subvariety and $K\subset Y$, then
$\PH K\subset Y$.  That is, $\PH K$ is contained in the Zariski closure of
$K$.   Furthermore, if $Y\subset \bbp^n$ is a projective manifold and
$\l=\co(1)\bigr|_Y$, the $\l$-projective hull of $K\subset Y$, defined
as in (2.1) with  $\co(1)$ replaced by $\l$,  agrees with $\PH K$.  
The same is true of $\lambda^k$  for any $k$.  
\medskip

\item {(2)}  If $K$ is contained in an affine open subset $\Omega\subset \bbp^n$,
then $(\PH K)_{{\rm poly}, \Omega}\subseteq \PH K$. 
\medskip

\item {(3)}   If $K=\partial C$, where $C$ is a holomorphic curve with boundary
in $\bbp^n$, then $C\subseteq \PH K$.
\medskip

\item {(4)}  ${\PH K}^- - K$ satisfies the {\sl maximum modulus principle}  for
holomorphic functions on open subsets of $\bbp^n-K$. (${\PH K}^-$ denotes the closure of ${\PH K}$.)
\medskip

\item {(5)}  ${\PH K}^- - K$ is 1-pseudoconcave in the sense of [DL]. In particular, for
any open subset $U\subset \bbp^n-K$, if the Hausdorff 2-measure 
$\ch^2(\PH K\cap U)<\infty$, then $\PH K \cap U$ is a {\sl complex analytic
subvariety} of dimension 1 in $U$.
\medskip

\item {(6)}  If $K$ is a real analytic curve, then the Hausdorff dimension of $\PH K $ is 2.
\medskip
\medskip

\item {(7)}     $\PH K$ is pluripolar if and only if $K$ is pluripolar.
\footnote{$^1$}{A set $K$ is pluripolar if  it is locally contained in the $-\infty$ set of 
a plurisubharmonic function.}
 Thus there exist
smooth closed curves   $\g \subset \bbp^2$ with $\PH\g=\bbp^2$.
However, real analytic curves are always pluripolar.
\medskip

The projective hull has simple characterizations in both affine and
homogeneous coordinates.  For example, if $\pi:\bbc^{n+1}-\{0\}\arr \bbp^n$
is the standard projection, and we set $S(K)\equiv \pi^{-1}(K)\cap
S^{1n+1}$, then 
$$
\PH K \ =\ \pi\biggl\{\PH{S(K)}_{{\rm poly}}-\{0\}\biggr\}.
$$
where $\PH{S(K)}_{{\rm poly}}$ is the polynomial hull of $S(K)$ in
$\bbc^{n+1}$. 

There is also a {\sl best constant function} $C:\PH K\to
\bbr^+$ defined at $x$ to be the least $C=C(x)$ for which the defining 
property (2.1) holds. For $x\in \PH K$, the set 
$\pi^{-1}(x)\cap \PH{S(K)}_{{\rm poly}}$ is a disk of radius
$\rho(x)=1/C(x)$. One deduces that $\PH K$ is compact if $C$ is bounded.

\Def{2.2} A compact subset $K\subset \bbp^n$ is called {\sl stable} if 
$C:\PH K \to \bbr^+$ is bounded.
\medskip

Combining a classical argument of E. Bishop with  (5) and (6) above  gives the following.

\Theorem{2.3. [HLW]} {\sl Let $\g\subset \bbp^n$ be a compact stable real analytic curve (not
necessarily connected).  Then $\PH \g -\g$ is a one-dimensional complex analytic
subvariety of $\bbp^n-\g$.}
\medskip

As mentioned in the introduction, there is much evidence for the conjecture that
any compact real analytic curve in $\bbp^n$ is stable, and therefore the stability hypothesis 
could be removed from Theorem 2.3. 
Interestingly,  the conclusion of Theorem 2.3 fails to hold in general for smooth curves (see
[HL$_5$, \S 4]), but may hold for smooth curves which are pluripolar 
 or, say, quasi-analytic.

\Remark{2.4} The close parallel between polynomial and projective hulls is 
signaled by the existence of a {\sl projective Gelfand transformation},
whose relation to the classical Gelfand transform is analogous to the
relation between Proj of a graded ring and Spec of a ring in modern
algebraic geometry. To any Banach graded algebra $A_*=\bigoplus_{d\geq0} A_d$
(a normed graded algebra which is a direct sum of Banach
spaces) one can associate  a topological space $X_{A_*}$ and a hermitian line
bundle $\l_{A_*}\to X_{A_*}$ with the property that $A_*$ embeds as a
closed subalgebra 
$$
 A_*\ \subseteq \ \bigoplus_{d\geq0} \Gamma_{{\rm cont}}(X_{A_*}, \l_{A_*}^d) 
$$
 of the algebra of continuous sections of powers of $\l$
with the sup-norm. When $K\subset \bbp^n$ is a compact subset and 
$A_d = H^0(\bbp^n, \co(d))\bigr|_K$ with the sup-norm on $K$, there is a
natual homeomorphism
$$
\PH K \ \cong \ X_{A_*}\and \co(1) \ \cong \ \l_{A_*}.
$$
This parallels the affine case where, for a compact subset $K\subset
\bbc^n$, the Gelfand spectrum of the  closure of the polynomials on
$K$ in the sup-norm corresponds to the polynomial hull of $K$.

Furthermore, when $A_*$ is finitely generated, the space $X_{A_*}$ 
can be realized, essentially uniquely, by a subset  $X_{A_*}\subset \bbp^n$
with $\l_{A_*}\cong \co(1)$ and $$\PH X_{A_*}\ =\ X_{A_*}.$$
This parallels the classical correspondence between finitely generated Banach
algebras and polynomially convex subsets of $\bbc^n$. Details are given in
[HL$_5$].
\medskip

\vfill\eject

\centerline{\bf \PL. Projective Linking and Projective Winding Numbers}
\medskip

In this section we introduce the notion of projective linking numbers and 
projective winding numbers for oriented curves in $\bbp^n$.

Let $M=M^{2p-1}\subset \bbp^n$ be a compact oriented submanifold of dimension
$2p-1$, and recall the {\sl intersection pairing}
$$
\bullet :H_{2p}(\bbp^n, M;\bbz) \,\times \, H_{2(n-p)}(\bbp^n-M;\bbz)\  \arr\
\bbz 
\eqno(\PL.1)
$$
which under Alexander Duality corresponds to the Kronecker pairing:
$$
\kappa :H^{2(n-p)}(\bbp^n- M;\bbz) \,\times \, H_{2(n-p)}(\bbp^n-M;\bbz)\ 
\arr\ \bbz.
$$
When homology  classes are represented by cycles which intersect
transversally in regular points, the map (\PL.1) is given by the usual algebraic intersection
number.

\Def{\PL.1} Fix $M$ as above and let $Z\subset \bbp^n-M$ be an algebraic
subvariety of dimension $n-p$.  Then the {\sl projective linking number}
of $M$ and $Z$ is defined as follows.  Choose a $2p$-chain $N$  in $\bbp^n$
with  $dN=M$ and set: 
$$
\Link_{\bbp}(M, Z)\ \equiv \ N\bullet Z - \deg(Z) \int_N \o^{p}
\eqno(\PL.2)
$$
where $\o$ is the standard K\"ahler form on $\bbp^n$. 
 This definition
extends by linearity to algebraic $(n-p)$-cycles $Z=\sum_j n_j Z_j$
supported in $\bbp^n-M$. 

\Lemma{\PL.2} {\sl The linking number $\Link_{\bbp}(M, Z)$ is independent of
the choice of the cobounding chain $N$.}

\pf  Let $N'$ be another choice and set $W=N-N'$.  Then $d W=0$ and 
$$
\qquad\qquad
W\bullet Z- \deg(Z)\int_{W}\o^{p} \ = \ \deg(W)\cdot\deg(Z) -
\deg(Z)\cdot\deg(W)\ =\ 0.   
\qquad \qquad\vrule width5pt height5pt depth0pt
$$ 

For now we shall be concerned with the case $p=1$. Here there is a naturally
related  notion  of projective winding number defined as follows.

\Def{\PL.3} Suppose $\G\subset \bbp^n$ be a smooth closed oriented curve and
let $\sigma$ be a holomorphic section of $\co(\ell)$ over $\bbp^n$ which does
not vanish on $\G$. Then the {\sl projective winding number} of $\sigma$ on
$\G$ is defined to be $$
\Wind_{\bbp}(\G, \s) \ \equiv\ \int_{\G} d^C \log \|\sigma\|
\eqno(\PL.3)
$$

\Prop{\PL.4}  {\sl  Let $\G$ and $\sigma$ be as in Definition \PL.3, and let
$Z$ be the divisor of $\sigma$.  Then }
$$
\Wind_{\bbp}(\G, \s) \ =\ \Link_{\bbp}(\G, Z) \ 
$$

\pf
We recall the fundamental formula
$$
dd^C \log\|\sigma\|\ =\ Z- \ell\o
\eqno(\PL.4)
$$
which follows by writing $\|\sigma\|= \rho |\sigma|$ in a local holomorphic
trivialization of $\co(\ell)$ and applying the Chern and Poincar\'e-Lelong
formulas: $dd^C \log \rho = -\ell\o$ and $dd^C\log|\s|= \Div(\s) = Z$.
We now write $\G=\partial N$ for a rectifiable 2-current $N$ in $\bbp^n$
and note that
$$
\int_{\G} d^C \log \|\s\|\ =\ \int_Ndd^C\log \|\s\|\ =\ 
N\bullet Z- \deg(Z)\int_N \o
$$
by formula (\PL.4).\qed
\medskip
It will be useful in what follows to normalize these linking numbers.

\Def{\PL.5} For $\Gamma$ and $Z$ as above, we define the {\sl reduced
projective linking number} to be
$$
\widetilde{\Link}_{\bbp}(\G, Z) \ \equiv\ {1\over {\deg(Z)}}
 \Link_{\bbp}(\G, Z)
$$
 and the {\sl reduced projective winding number} to be
$\widetilde{\Wind}_{\bbp}(\G, Z)\equiv {\Wind}_{\bbp}(\G, Z)/\deg(Z)$.
\smallskip
Note that if $Z=\Div(\sigma)$ for $\s \in H^0(\bbp^n,\co(\ell))$, then by
Proposition \PL.4 we have
$$
\widetilde{\Link}_{\bbp}(\G, Z) \ =\ \int_{\G} d^C \log \|\s\|^{1\over
{\ell}} 
\eqno(\PL.5)
$$

\Remark{\PL.6} {\bf (The relation to affine linking numbers.)} Note that
a smooth curve $\G\subset\bbp^n$ does not meet the generic hyperplane
$\bbp^{n-1}$ and is therefore contained in the affine chart 
$\bbc^n=\bbp^n-\bbp^{n-1}$.  Let $N$ be an integral 2-chain with support in
$\bbc^n$ and $dN=\G$.  Then 
$$
\Link_{\bbp}(\G, \bbp^{n-1}) \ =\ N\bullet \bbp^{n-1}
- \int_N\o
\ =\ - \int_N\o,
\eqno(\PL.6)$$
and for any divisor $Z$ of degree $\ell$ which does not meet $\G$ we have
$$
\Link_{\bbp}(\G, Z) -\Link_{\bbp}(\G, \ell\bbp^{n-1}) \ =\ N\bullet Z
 \ =\ \Link_{\bbc^n}(\G, Z),
$$
the classical linking number of $\G$ and $Z$.

\Remark{\PL.7} {\bf (The relation to sparks and differential characters.)}
The deRham-Federer approach to differential characters (cf. [GS], [Ha],
[HLZ]) is built on objects called {\sl sparks}. These are generalized
differential forms (or currents) $\a$  which satisfy the equation
$$
d\a \ =\ R-\phi
$$
where $R$ is rectifiable and $\phi$ is smooth.
By (\PL.4) the generalized 1-form $d^C\log\|\s\|$ is such a creature.
Its associated Cheeger-Simons character $[\a]:  Z_1 (\bbp^n)\arr\bbr/\bbz$
on smooth 1-cycles is defined by the formula 
$$
[\a] (\Gamma)\ \equiv\ \int_N \ell\omega\ \ \ ({\rm mod}\ \ \bbz) 
$$
where $N$ is a 2-chain with  $dN=\G$.  Thus we see that the ``correction term''
in the projective linking number is the value of a secondary invariant
related to differential characters.

\Remark{\PL.8} Note that by
(\PL.4) the function $u\equiv \log\|\s\|^{1\over{\ell}}$ which appears in
(\PL.5) has the property that  $$
dd^C u +\o\ ={1\over{\ell }} Z\ \geq\ 0
\eqno(\PL.7)$$
by the fundamental equation (\PL.4).  This leads us naturally to the next
section.

\vskip .3in

\centerline{\bf \PSH. Quasi-plurisubharmonic Functions}
\medskip
  
The following concept, due to Demailly [D$_1$] and developed systematically
by  Guedj and Zeriahi [GZ], is central to the study of projective hulls,
and it is intimately related to projective linking numbers.

\Def{\PSH.1}  An upper semi-continuous function $u:X \arr [-\infty, \infty)$ on
a K\"ahler manifold $(X,\omega)$ is {\sl \omp} (or {\sl $\omega$-quasi-plurisubharmonic})
  if $u\not\equiv -\infty $  and
$$
d d^C u + \o\ \geq 0\qquad  {\rm on}\ \ X.
\eqno(\PSH.1)
$$
The convex set of all such functions on $X$  will be denoted by  $\psh(X)$.
\medskip

These functions enjoy many of the properties of classical plurisubharmonic
functions and play an important role in understanding various capacities in
projective space (cf. [GZ]). One of the appealing geometric properties of
this class is the following. Suppose $\o$ is the curvature form of a
holomorphic  line bundle $\l\to X$ with hermitian metric $g$. Then a smooth
function $u:X\to \bbr$ is \omp iff the hermitian metric $e^u g$ has
curvature $\geq 0$ on $X$.
  
The \omp functions are directly relevant to projective hulls, 
as the next result shows (cf. [GZ, Pf. of Thm. 4.2] and Theorem 5.3 below).

\Theorem{\PSH.2}  {\sl Let $\o$ denote the standard K\"ahler form on $\bbp^n$.
Then the projective hull $\PH K$ of a compact subset $K\subset \bbp^n$, is
exactly  the subset of points $x\in\bbp^n$ for which there exists a constant
$\L =\L(x)$ with}
$$
 u(x) \ \leq \sup_K u + \L\qquad {\rm for\ all}  \ \ u \in \psh(X)
\eqno(\PSH.2)$$

\medskip
This    enables one to generalize the notion of projective hull from
projective algebraic manifolds to general K\"ahler manifolds.

Note that the least constant $\L(x)$ for which (\PSH.2) holds is exactly
$
\L(x)=\log C(x)
$
where $C(x)$ is the best constant function discussed in \S2.

By considering functions of the form $u=\log\{\|\sigma\|^{1\over d}\}$
for  sections $\sigma\in H^0(\bbp^n, \co(d))$, one immediatly sees the
necessity of the condition in 
Theorem \PSH.2. Sufficiency follows from the fact that such
functions are the extreme points of the cone $\psh(\bbp^n)$.

One importance of  \omp functions is that they enable us to establish
 Poisson-Jensen measures for points in $\PH K$, [HL$_5$, Thm. 11.1]. 
\vskip .3in

\centerline{\bf \PHC. Boundaries of Positive Holomorphic Chains.}
\medskip

Let $\G$ be a smooth oriented closed curve in $\bbp^n$.  We recall that
(even if $\G$ is only class $C^1$) any irreducible complex analytic
subvariety $V\subset \bbp^n-\G$
of dimension one has finite Hausdorff 2-measure and defines
a current $[V]$ of dimension 2 in $\bbp^n$ by integration on the
canonically oriented manifold of regular points.  This current satisfies
$d[V]=0$ in $\bbp-\G$. (See [H] for example.)

\Def{\PHC.1} By a {\sl positive holomorphic 1-chain with boundary $\G$} we
mean a finite sum $T=\sum_k n_k [V_k]$ where each $n_k\in\bbz^+$ and
each $V_k\subset \bbp^n-\G$ is an irreducible subvariety of dimension 1, so 
that
$$
dT\ =\ \G\qquad\ \ {\rm (as\ currents\ on }\ \  \bbp^n)
$$

We shall be interested in conditions on $\G$ which are necessary and
sufficient for it to be such a boundary.

\Prop{\PHC.2}  {\sl Let $\G\subset \bbp^n$  be a smooth oriented closed
curve (not necessarily connected) with a positive integer multiplicity on each component. Then the following are equivalent.
\medskip

(1) \qquad $\widetilde{\Link}(\G, Z)\ \geq \ -\L$ \ \ for all algebraic
hypersurfaces $Z\subset  \bbp^n-\G$.

\medskip

(2) \qquad $\widetilde{\Wind}(\G, \s)\ \geq \ -\L $ \ \ for all
$\s\in H^0(\bbp^n, \co(\ell)),  \ell>0$, \  \  with no zeros on $\G$.

\medskip

(3) \qquad $\int_\G d^C u\ \geq\ - \L$\qquad\ \ \  for all $u \in
\psh(\bbp^n)$.
}

\pf
The equivalence of (1) and (2) is follows immediately from Proposition 
\PL.4. That (3) $\Rightarrow$ (1) follows from (\PL.5) and (\PL.7).
To see that (1) $\Rightarrow$ (3) we use the following non-trivial fact due
essentially to Demailly [D$_2$], [G].

\Prop{\PHC.3} {\sl The functions of the form $\log\|\s\|^{1\over{\ell}}$
for $\s\in  H^0(\bbp^n, \co{}({\ell}))$, $\ell  > 0$,   are weakly dense in
$\psh(\bbp^n)$ modulo the constant functions.}  

\pf
Fix $u\in \psh(\bbp^n)$ and consider the positive (1,1)-current 
$T\equiv d d^C u +\omega$.
By [D$_2$], [G, Thm. 0.1],  there exist sequences 
$\s_j\in H^0(\bbp^n, \co{}(\ell_j))$ and $N_j\in \bbr^+$, $j=1,2,3,...$ such
that  
$$
T\ =\ \lim_{j\to\infty}{1\over N_j} \Div(\s_j)\ =\ 
\lim_{j\to\infty}{\ell_j\over N_j}\bigg\{ d d^C\log\|\s_j\|^{1\over \ell_j}
+\omega\bigg\}.
\eqno{(\PHC.1)}
$$
Since $T(\omega^{n-1}) = \int_{\bbp^n}(dd^Cu\wedge \omega^{n-1}+ \o^n)
 = \int_{\bbp^n} \o^n  =1 =\lim_{j\to \infty} \ell_j/N_j$ by (\PHC.1), 
 we may assume that $N_j=\ell_j$ for all $j$.
Therefore, setting $u_j = d d^C\log\|\s_j\|^{1\over \ell_j}$, 
we have   by (\PHC.1) that $\lim_jdd^C u_j  = dd^C u$  and, in particular, 
$\lim_j \D u_j=\D u$.
Renormalizing each $u_j$ by an additive constant, we may also assume that $\lim_j\int u_j = \int u$.
 
Recall the formula $Id =  H + G\circ \D$ on $C^\infty(X)$ where $H$ is harmonic
projection and $G$ is the Green's function. This decomposition carries over to distributions 
$\cd'(X)$ by adjoint. It follows that $u_j-\int u_j = G(\D u_j) \to G(\D u) = u-\int u$ and so 
$u_j\to u$ as claimed.
\qed

\Theorem{\PHC.4} {\sl Let $\G$ be as above. Suppose $\G = dT$  where $T$ is
a   positive holomorphic chain with mass $\M(T) \leq \L$. Then $\G$ satisfies
the equivalent conditions (1), (2), (3) of Proposition \PHC.2. 
}

\pf Suppose $u\in \psh(\bbp^n)$.  Since $dd^C u +\o \geq 0$, we have  
$$
\int_\G d^C u \ =\ \int_{dT} d^C u \ =\ T(dd^Cu)\ =\ T(dd^Cu+\o) -T(\o)
\ \geq\ -T(\o)\ =\ -\M(T)\ \geq-\L
$$
as asserted. \qed


\vskip .3in

\centerline{\bf \PAW. The Projective Alexander-Wermer Theorem for Curves.} 
\medskip
  
We now examine the question of whether the equivalent necessary conditions,
given in Theorem \PHC.4  are in fact sufficient.

\Theorem{\PAW.1} {\sl  Let  $\G\subset \bbp^n$  be an  embedded, oriented,   real analytic  closed
curve, not necessarily connected and with  a positive   integer multiplicity
on each component.  Suppose the underlying curve $\supp(\G)$ is stable.
  If $\G$ satisfies the
equivalent conditions (1), (2), (3) of Proposition \PHC.2, then 
there exists a positive holomorphic 1-chain $T$ in $\bbp^n$ with $\M(T)\leq
\Lambda$ such that }
$$
dT \ =\ \G.
$$

\noindent
{\AAA N\BBB OTE.} When it is clear from context, we   also use  $\G$ to denote
the underlying curve supp$\,\G$.

\pf Fix a compact set $K\subset \bbp^n$ and consider the sets
$$\eqalign{
C_K\ &\equiv\ \{ dd^CT \ : \ T\in\cp_{1,1},\ \  \M(T)\leq 1\ \ {\rm and}\ \ 
\supp\, T\subset K\}  
\cr
\cs_K \ &\equiv \ \{u\in C^\infty(\bbp^n)\ :\ dd^C u + \o \geq 0\ \  {\rm
on}\  K\} }
$$
where $\cp_{1,1}\subset\ce_2(\bbp^n)$ denotes the convex cone of positive
currents of bidimension (1,1) on $\bbp^n$. Note that $C_K$ is a weakly
closed, convex subset which contains 0 (since the set of 
$T\in\cp_{1,1}$ with $\M(T)\leq 1$ and  $\supp\, T\subset K$ is weakly
compact). 

Suppose that   $\ck$ is a weakly closed convex subset containing 0 in a
topological vector space $V$. Then the {\sl polar} of $\ck$ is defined to
be the subset of the dual space given by
$$
\ck^0\ \equiv\ \{L\in V' : (L,v)\ \geq\ -1\ \ {\rm for }\ \ v\in\ck\}.
$$
Similarly given   a subset $\cl\subset V'$ we define
$$
\cl^0\ \equiv\ \{v\in V :  (L,v) \ \geq\ -1\ \ {\rm for }\ \ L\in\cl\}.
$$
The Bipolar Theorem [S] states that 
$$
(\ck^0)^0 \ =\ \ck.
$$

\Prop{\PAW.2}
$$
\cs_K\ =\ (\cc_K)^0
$$
 \pf Suppose that $u\in  C^\infty(\bbp^n)$ satisfies 
$$
u(dd^C T) \ =\ T(dd^C u) \ \geq\ -1
\eqno(\PAW.1)$$
for all $T\in\cp_{1,1}$ with $\M(T)\leq 1$ and  $\supp\, T\subset K$. 
Consider $T=\d_x \x$ where   $x\in K$ and $\x$ is a positive (1,1)-vector
 at $x$ with mass-norm $\|\x\|=1$. Then 
$$
T(dd^C u) \ =\ (dd^C u)(\x)
 \ =\ (dd^C u+\o)(\x) - 1\ \geq\ -1
$$
by (\PAW.1), and so $u \in \cs_K$. This proves that $(\cc_K)^0\subset
\cs_K$.

For the converse, let $u\in \cs_K$ and fix $T\in \cp_{1,1}$
with $\M(T)\leq 1$ and  $\supp\, T\subset K$. Then
$$
(dd^C T)(u)\ =\ T( dd^C u)\ =\ T( dd^C u+\o) -T(\o)\ \geq\ T(\o)
\ =\ -\M(T)\ \geq \ -1,
$$
and so $u\in (\cc_K)^0$. \qed
\medskip

As an immediate consequence we have the following. If we set
$$\eqalign{
\L\cdot\cs_K \ &\equiv\  \{dd^C T : T\in\cp_{1,1}, \M(T)\leq \L \ {\rm and}\
\supp\ T\subset K\}, \qquad{\rm then} 
\cr
(\L\cdot\cs_K)^0\ &=\ \{u\in C^{\infty}(\bbp^n)\ :\ dd^Cu + \o \geq -\L
\ \ {\rm on\ \ } K\}
}
\eqno(\PAW.2)$$

Recall from \S \PSH \ \ that for each $\L\in\bbr$ we have the compact set
$$
\PH{\G}(\L)\ =\ \{x\in \bbp^n\ :\ u(x) \leq \sup_{\G}\, u +\L \ \ 
\forall u \in \psh(\bbp^n)\} 
$$
and that $\PH{\G}=\bigcup_{\L}\PH{\G}(\L)$.  The following lemma is
established in [HL$_5$, 18.7].

\Lemma{\PAW.3}  {\sl  Let $u$ be a $C^\infty$ function which is defined and
\omp on a neighborhood of $\{\PH{\G}\}^-$. Fix $\Lambda\in
\bbr$.  Then there is a 
 $C^\infty$ function $\widetilde u$ which is defined and
\omp on all of $\bbp^n$ and agrees with $u$ on a
neighborhood of $\PH{\G}(\L)$. }
\medskip

Here $\{\PH{\G}\}^-$ denotes the closure of $\PH{\G}$. 
By our stability assumption, $\PH{\G}=\PH{\G}(\L)$ for some 
$\L$, and therefore  $\PH{\G}=\{\PH{\G}\}^-$.
 However, we shall keep the notation 
$\{\PH{\G}\}^-$, when appropriate, in order to prove the more general 
result mentioned in Remark \PAW.5 below.

The lemma above leads to the following.

\Prop{\PAW.4}  {\sl  If $\G$ satisfies condition (3) in Proposition \PHC.2,
that is, if 
$$
(-d^C \G)(u) \ =\ \int_{\G} d^C u\ \geq \ -\L
$$
for all $u\in\psh(\bbp^n)$, then there exists $T\in \cp_{1,1}$ with
$\M(T)\leq \L$ and $\supp T\subset \{\PH{\G}\}^-$, such that }
$$
dd^C T\ =\ -d^C \G
\eqno(\PAW.3)
$$
 \pf
If this were not true, then by compactness it would fail on some compact
neighborhood $K$ of $\{\PH{\G}\}^-$. By   (\PAW.2)  and the Bipolar
Theorem, we conclude that there is a smooth function $u$ which is
\omp on $K$ with $-d^C\G(u)<-\L$. Applying Lemma \PAW.3
contradicts the hypothesis. \qed
\medskip

Now let $T$  be the  current    given by Proposition \PAW.4.
Let $V$ denote the projective hull of $\G$ and recall from [HLW, Theorem 4.1]  that
$V$ has the following strong regularity.  There exists a Riemann surface $S$ with finitely
many components, a compact region $W \subset S$ with real analytic boundary, 
and a holomorphic  map $\rho:S\to \bbp^n$ which is generically injective and satisfies
  \medskip\noindent
  (1) \ \ $\rho(  W)\ =\   V$, 
  \medskip\noindent
  (2) \ \  $\rho$ is an embedding on a tubular neighborhood of $\partial   W$ in $S$,
  and 
  \medskip\noindent
  (3) \ \  $\rho(\partial   W)$ is a union of components of the support of $\G$.

\medskip

Let $\Sigma$ denote an $\e$-tubular neighborhood of $\partial W$ on $S$ (with respect to some 
analytic metric), with $\e$ chosen   so that 
 \medskip\noindent
  (4) \ \  $\rho$ is injective on $\Sigma$.
  \medskip\noindent
  (5) \ \  $\rho(\partial^+ \Sigma) \cap V\ =\ \emptyset$ where $\partial^+ \Sigma$ denotes the
  ``outer'' boundary of $\Sigma$, i.e. the union of 
  
  components of $\partial \Sigma$ not contained
  in $W$.
 \medskip
 \def\coeff{\nu}

 Write $\Sigma = \Sigma^+\cup \partial W\cup\Sigma ^-$ where $\Sigma^+$ denotes the union 
 of components of $\Sigma -\partial W$ which are not contained in $W$.  Then we have
 $d[\Sigma^+] = [\partial^+\Sigma] - [\partial W]$ (with standard orientations coming from the
 domains).  Hence we have
 $$
 d d^C  [\Sigma^+] \ =\ -d^C [\partial^+\Sigma] + d^C [\partial W].
 $$
Let $\coeff :\Sigma\to \bbz$ denote the locally constant, integer-valued function with
the property that
$$
\rho_*\left(\coeff  [\partial W]\right)\ =\ \Gamma.
$$
Then 
$$
 d d^C  \rho_*\left(\coeff  [\Sigma^+]\right)
  \ =\ -d^C \rho_*\left(\coeff  [\partial^+\Sigma]\right) + d^C\Gamma
   \ \equdef\ -d^C\Gamma^++d^C\Gamma.
\eqno(\PAW.4)
$$
We now define
$$
T^+\ =\ T + \rho_*\left(\coeff  [\Sigma^+]\right)
$$
and note that by (\PAW.3) we have
$$
dd^C T^+\ =\ -d^C \Gamma^+.
\eqno(\PAW.5)
$$
Now  in the open set $\bbp^n-  \G^+$,  the current $T^+$ is  a positive (1,1)-current which is supported in the analytic subvariety $V^+\equiv \rho(W^+)\equiv \rho(\Sigma^+\cup \overline W)$ and
satisfies  $dd^CT^+=0$.  We note that $\rho:W^+\to V^+$ is just the normalization of $V^+$.
It follows from [HL$_3$, Lemma 32] that there is a harmonic function $h:W\to\bbr^+$
so that 
$$
T^+ = \rho_*(h[W^+]) \qquad {\rm in}\ \   \bbp^n-\G^+.
$$
This function  is evidently constant ($= \coeff $) in $\Sigma^+$, and 
hence $h$ is constant on every component of $W^+$.  Thus $T= \rho_*(h[W])$
is a positive holomorphic chain,

Now since $\supp dT\subset \supp \G$ and $\supp \G$ is a 
regular curve, it follows from the Federer Flat Support Theorem
[F, 4.1.15] that $dT$ is of the form $dT=\sum_k c_k [\G_k]$
where the $c_k$'s are constants and the $\G_k$'s are the oriented 
 connected components of $\supp\G$. The current $\G$ itself has the form
 $\G=\sum_k \coeff _k [\G_k]$ for integers $\coeff_k$. Since $dd^C T = -d^C\Gamma$, 
 we have $d^C(dT-\G) = \sum (c_k-\coeff _k)d^C[\G_k]=0$.  It follow that $c_k=\coeff_k$ for all $k$, 
 that is, $dT=\G$ as claimed.\qed

\Remark{\PAW.5}  Theorem \PAW.1 remains true if one replaces the second statement
with the assumption that  the closed projective hull $\{\PH \G\}^-$ of
$\G$ is locally contained in a 1-dimensional complex
subvariety at each point of $\{\PH \G\}^- -\G$.  
\medskip

Theorem \PAW.1 can be restated in terms of the classical (affine) linking
numbers.

\Theorem {\PAW.6}  {\sl  Let $\G  \subset \bbc^n\subset \bbp^n$ be an oriented
closed curve as in Theorem \PAW.1.
Then there exists a constant $\Lambda_0$ so that the classical linking
number 
$$
\Link_{\bbc^n}(\G, Z) \ \geq\ -\Lambda_0\ \deg Z
$$
for all algebraic hypersurfaces in $\bbc^n-\G$  if and only if 
$\G$ is the boundary of a positive holomorphic 1-chain $T$ in $\bbp^n
= \bbc^n\cup \bbp^{n-1}$  with (projective) mass }
$$
\M(T)\ \leq \ \Lambda_0 + {1\over2}\int_{\G} d^C\log ({1+\|z\|^2})
$$
\pf
This follows from Remark \PL.6 and the fact that the projective K\"ahler
form $\omega ={1\over2}dd^C\log(1+\|z\|^2)$ in $\bbc^n$.\qed

\Note{\PAW.7} The case $\Lambda_0=0$ corresponds to  the theorem of Alexander
and Wermer [AW$_2$]. A proof of the full Alexander-Wermer Theorem
(for $C^1$-curves with no stability assumption) in the spirit of the arguments above
is given in [HL$_6$].
\medskip

Note that the current of least mass among those provided by Theorem \PAW.1 is uniquely
determined by $\G$. All others differ from this one by adding a positive algebraic cycle. 

\Cor{\PAW.8} {\sl  Let $\g$ be as in Theorem \PAW.1 and suppose  $T$ is a
 positive holomorphic 1-chain with $dT=\g$.  Then $T$ is the unique holomorphic chain 
 of least mass with  $dT=\g$ if and only if  
 $$
 \inf_Z \left\{ {T\bullet Z \over \deg Z} \right\} \ =\ 0
 $$
where the infimum is taken over all algebraic hypersurfaces in the complement $\bbp^n-\g$.}
\medskip

\pf Suppose $\inf_Z \{T\bullet Z/\deg Z\} = c>0$.  Then $\widetilde{\Link}_{\bbp}(\G, Z) = T\bullet Z/\deg(Z)-T(\o) \geq c-T(\o)$ for all positive divisors in $\bbp^n-\G$.  Hence, by Theorem \PAW.1 there exists a
positive holomorphic chain $T'$ with $dT'=\G$ and $M(T')\leq M(T)-c$. 

On the other hand suppose that $T$ is not the unique positive holomorphic chain $T_0$ of least mass.
Then $T=T_0+ W$ where $W$ is a positive algebraic cycle, and one has
$ T\bullet Z/\deg Z=  T_0\bullet Z/\deg Z + W\bullet Z/\deg Z \geq  W\bullet Z/\deg Z=\deg W$.\qed

\vfill\eject


\centerline{\bf \PM. Theorems for General Projective Manifolds.}  \medskip
  
The results established above generalize from $\bbp^n$ to any projective
manifold.  Let $X$ be a compact complex manifold with a positive 
holomorphic line bundle $\lambda$.
Fix a hermitian metric on $\lambda$ with curvature form $\omega>0$, and
give $X$ the K\"ahler metric associated to $\omega$.
Let $\G$ be a closed curve with integral weights as in Theorem \PAW.1,
and assume $[\G]=0$ in $H_1(X;\,\bbz)$.

\Def{\PM.1}  Let $Z=\Div(\s)$ be the divisor of a holomorphic section
$\s \in H^0(X, \co{}(\lambda^\ell))$ for some $\ell\geq 1$. If $Z$
does not meet $\G$, we can define the  {\bf linking number} and the {\bf
reduced linking number} by
$$
\Link_{\lambda}(\G,Z)\ \equiv \ N\bullet Z-\ell\int_N\omega
\and
\widetilde{\Link}_{\lambda}(\G,Z)\ \equiv \ {1\over \ell}\ \Link(\G,Z)
$$
respectively, where $N$ is any 2-chain in $Z$ with $dN=\G$ and where
the intersection pairing $\bullet$ is defined as in (\PL.1) with $\bbp^n$
replaced by $X$.

\medskip

To see that this is well-defined suppose that $N'$ is another 2-chain
with $dN'=\G$. Then
 $(N-N')\bullet Z - \ell\int_{N-N'}\omega = (N-N')\bullet
(Z-\ell[\omega])=0$ because $Z-\ell\omega$ is cohomologous to zero in $X$.

\Def{\PM.2}  The {\bf  winding number} of a section 
$\s \in H^0(X, \co{}(\lambda^\ell))$ with $\|\s\|>0$ on $\G$,  is
defined to be
$$
\Wind_{\lambda}(\G, \s)\ \equiv\ \int_\G d^C \log\|\s\|.
$$
The {\bf reduced winding number} is
$$
\widetilde{\Wind}_{\lambda}(\G,\s)\ \equiv\ {1\over\ell  } \ 
\Wind_{\lambda}(\G, \s) \ =\ \int_\G d^C\log\|\s\|^{1\over\ell}
$$

From the Poincar\'e-Lelong equation
$$
d d^C \log\|\s\|\ =\ \Div(\s) -\ell \omega
\eqno(\PM.1)
$$
we see that
$$
\Wind_{\lambda}(\G, \s)\ =\ \Link_{\lambda}(\G, \Div(\s)).
\eqno(\PM.2)
$$
From (\PM.1) we also see that $\log\|\s\|^{1\over\ell}$ belongs to the
class $\psh(X)$ of \omp functions
on $X$ defined in  \PSH.1.

\Prop{\PM.3}  {\sl The following are equivalent:

\medskip

(a) \   $\widetilde{\Link}_{\lambda}(\G,Z)\ \geq \ -\Lambda$ \ 
for all divisors $Z$ of holomorphic sections of $\lambda^\ell$, and all $\ell>0$.

\medskip

(b) \   $\widetilde{\Wind}_{\lambda}(\G,\s)\ \geq \ -\Lambda$ \ 
for all holomorphic sections $\s$ of $\lambda^\ell$, and all
 $\ell>0$.
\medskip

(c) \ \  $\int_{\G} d^Cu\ \ \geq \ -\Lambda$ \ \ for all $u\in
\psh(X)$.

}
\pf  That (a) $\Leftrightarrow$ (b) $\Rightarrow$ (c) is clear.
To see that (c) $\Rightarrow$ (b) we use results of Demailly. Fix 
$u\in \psh(X)$ and consider the positive closed (1,1)-current
 $T\equiv d d^C u +\omega$. Note that $[T] = [\omega] = c_1(\lambda)
\in H^2(X;\bbz)$.  It follows from [D$_2$] that $T$ is the weak limit 
$$
T=\lim_{j\to \infty}{1\over N_j}\Div(\s_j)
$$
where $\s_j\in H^0(X,\co{}(\lambda^{\ell_j}))$ and $N_j>0$.
We can normalize this sequence by scalars so that
$\M({1\over N_j}\Div(\s_j))=\M(T)$ for all $j$. 
Set $\Omega= {1\over (n-1)!}\o^{n-1}$. Then $\M(T) =
\M({1\over N_j}\Div(\s_j)) = ({1\over N_j}\Div(\s_j), \Omega)
= {1\over N_j}(\ell_j\omega, \Omega)
= {\ell_j\over N_j}(\omega, \Omega)
={\ell_j\over N_j}\M(T)$.
 Therefore, $\ell_j=N_j$ for all
$j$.  Since $\Div(\s_j) = dd^C\log\|\s_j\| + \ell_j\omega$, we conclude that
$$
d d^C u = T-\omega = \lim_{j\to\infty}\bigg\{d d^C
\log\|\s_j\|^{1\over \ell_j} + \omega\bigg\}-\omega
\ =\  
 \lim_{j\to\infty}d d^C
\log\|\s_j\|^{1\over \ell_j}.$$
The remainder of the argument replicates the one  given for Proposition \PHC.3.
\qed

\medskip

For a compact subset $K\subset X$ the authors introduced the notion of
the $\lambda$-hull $\widehat{K}_{\lambda}$ of $K$ and showed that 
 $\widehat{K}_{\lambda}=\PH K$ for any embedding $X\hookrightarrow \bbp^N$ 
by sections of some power of $\lambda$.  We shall say that $K$ is {\bf $\lambda$-stable}
if $\widehat{K}_{\lambda}$  is compact.

\Theorem{\PM.4}  {\sl   Let $\Gamma\ =\ \sum_{
a=1}^M m_\a \G_\a$ be an embedded, oriented,  real analytic closed curve with 
integer multiplicities in $X$, and assume $\G$ is $\lambda$-stable.  Then there exists a positive holomorphic 1-chain
$T$ in $X$ with $dT=\G$ and $\M(T)\leq \Lambda$ if and only if any of the
equivalent conditions of Proposition \PM.3 is satisfied.}

\pf
If $dT=\G$ and $\M(T)\leq \Lambda$, then (a) follows as in the proof of
\PHC.4 above.

For the converse we recall from [HL$_5$,(4.4)] that the hull
$\widehat{K}_{\lambda}$ of a compact subset $K\subset X$  consists exactly
of the points $x$ where the extremal function
$$
\Lambda_K(x) \ \equiv\ \sup\{u(x) : u\in\psh(X)
\ {\rm and\ }\ u\leq0\ {\rm on }\ X\}
$$
is finite.  This enables one to directly carry through the arguments for
Theorem \PAW.1 in this case.
\qed

\medskip

\vskip .3in

\centerline{\bf \MH.  Relative holomorphic cycles.}
\medskip

Our main theorem \PAW.1 has a nice reinterpretation in terms of 
the Alexander-Lefschetz duality pairing discussed in section \PL.

  \Theorem{8.1} {\sl Let  $\gamma\subset \bbp^n$  be  a finite disjoint union of
  real analytic curves with $\gamma$ stable. Then a class
   $\tau\in H_2(\bbp^n,\gamma; \bbz)$ is represented by a positive holomorphic chain
   with boundary  on $\gamma$ 
   if and only if 
   $$
   \tau\bullet u\ \geq\ 0 
  \eqno{(8.1)}
   $$
  for all 
   $u\in H_{2n-2}( \bbp^n-\gamma;\bbz)$  represented by  positive algebraic hypersurfaces in $\bbp^n-\gamma$.}

  \pf  The implication $\Rightarrow$ is clear from the positivity of complex intersections.
  For the converse,   consider the short exact sequence
$$
0\ \arr \ H_2(\bbp^n; \bbz)\ \arr \ H_2(\bbp^n, \gamma; \bbz)\ \Arr{$\delta$} \ H_1(\gamma; \bbz)\ \arr \ 0.
$$
Note that 
$$
\delta \tau\ =\  \sum_{k=1}^\ell m_k {\oa \gamma}_k \ \equiv\ \G
$$
  where ${\oa \gamma}_1,..., {\oa \gamma}_\ell$ are the connected components of $\gamma$
  with a chosen orientation and the $n_k$'s are integers which we can assume to be positive.
   We are assuming that (8.1) holds for any positive algebraic class $u$.
    If $\G=0$, the desired conclusion is immediate, so we assume that $\G\neq 0$.
Now let $Z\subset \bbp^n-\gamma$ be any positive  algebraic hypersurface, and note that  
$$
0\ \leq \ {\tau\bullet Z  \over \deg Z} \ =\ 
{\tau\bullet Z  \over \deg Z}  -\tau(\o)+ \tau(\o) \ =\ \Link_{\bbp}(\G,Z) + \tau(\o).
$$ 
 Therefore, by Theorem \PAW.1, there exists a positive holomorphic 1-chain $T$
 with  $dT=\G$ and $ M(T) \leq \tau(\o)$.  Note that $\delta([T]-\tau)=0$ and $([T]-\tau)(\o)\geq0$.
Hence, $[T]-\tau$ is a positive class in $H_2(\bbp^n; \bbz)$ and is represented by a positive
 algebraic 1-cycle, say $W$.  Therefore, $\tau = [T+W]$ is represented by a positive
 holomorphic chain as claimed.\qed
 
\medskip
This result leads to a  nice duality between the cone in    
$H_2(\bbp^n, \gamma; \bbz)$ of  those classes which are 
represented by positive holomorphic chains, and the cone in
 $H_{2n-2}( \bbp^n-\gamma;\bbz)$ of classes 
represented by positive algebraic hypersurfaces.  
This duality, in fact,    extends to  more general projective
manifolds.  All this is discussed in detail  
in [HL$_8$]

 \vskip .3in

\centerline{\bf References}

\vskip .2in


 \noindent
[AW$_1$]  H. Alexander and J. Wermer,  Several Complex
Variables and Banach Algebras, Springer-Verlag,  New York, 1998. 

 \smallskip

\noindent
[AW$_2$]  H. Alexander and J. Wermer, {\sl Linking numbers
and boundaries of varieties}, Ann. of Math.
{\bf 151} (2000),   125-150.

 \smallskip

\noindent
[D$_1$]  J.-P. Demailly, {\sl Estimations $L^2$ pour
l'op\'erateur $\overline{\partial}$ d'un fibr\'e vectoriel holomorphe
semi-positif au-dessus d'une vari\'et\'e k\"ahl\'erienne compl\`ete},
   Ann. Sci. E. N. S.  
 {\bf 15} no. 4  (1982),  457-511.

 \smallskip

\noindent
[D$_2$]  J.-P. Demailly, {\sl Courants positifs extr\^emaux et
conjecture de Hodge},
 Inventiones Math.
 {\bf 69} no. 4  (1982),  347-374.

 \smallskip

\noindent
 [DL]   T.-C. Dinh and M. Lawrence,  {\sl
 Polynomial hulls and positive currents},
  Ann. Fac. Sci de Toulouse {\bf 12} (2003),   317-334.

 \smallskip

 \noindent
[Do]  P. Dolbeault,   {\sl On holomorphic chains with given 
boundary in $\bbc\bbp^n$},    Springer  Lecture Notes in Math., {\bf 
1089}  (1983), 1135-1140.

\smallskip

\noindent
[DH$_1$]   P. Dolbeault and G. Henkin ,  
{\sl  Surfaces de Riemann de bord donn\'e dans $\bbc\bbp^n$  },  pp.
163-187 in  `` Contributions to Complex Analysis and Analytic Geometry''  ,
  Aspects of Math.  Vieweg{\bf 26 },   1994.

\smallskip

\noindent
[DH$_2$]   P. Dolbeault and G. Henkin,   
{\sl Cha\^ines holomorphes de bord donn\'e dans $\bbc\bbp^n$},    
Bull.  Soc. Math. de France,
{\bf  125}  (1997), 383-445.

\smallskip

 \noindent
[DS] J. Duval and N. Sibony {\sl
Polynomial convexity, rational convexity and currents},
  Duke Math. J. {\bf 79}  (1995),     487-513.

 \smallskip

\noindent
[F]   H. Federer, Geometric Measure  Theory,
 Springer--Verlag, New York, 1969.

 \smallskip

\noindent
[GS] H. Gillet and C. Soul\'e,  {\sl
Arithmetic chow groups and differential characters},
pp. 30-68  in ``Algebraic K-theory;  Connections with Geometry and
Topology'',     Jardine and Snaith (eds.), Kluwer Academic Publishers,
  1989.
\smallskip

\noindent
[G]  V. Guedj,     
{\sl    Approximation of currents on complex manifolds},    
Math. Ann. {\bf 313} (1999), 437-474.

\smallskip

\noindent
[GZ]  V. Guedj and A. Zeriahi,     
{\sl    Intrinsic capacities on compact K\"ahler manifolds},    
Preprint Univ. de Toulouse , 2003

\smallskip

\noindent
[Ha]   B. Harris, {\sl
Differential characters and the Abel-Jacobi map},
pp. 69-86 in ``Algebraic K-theory;  Connections with Geometry and
Topology'', Jardine and Snaith (eds.), Kluwer Academic Publishers,
1989.
\smallskip

\noindent
[H]  F.R. Harvey,
Holomorphic chains and their boundaries, pp. 309-382 in ``Several Complex
Variables, Proc. of Symposia in Pure Mathematics XXX Part 1'', 
A.M.S., Providence, RI, 1977.

 \noindent 
\noindent
[HL$_1$] F. R. Harvey and H. B. Lawson, Jr, {\sl On boundaries of complex
analytic varieties, I}, Annals of Mathematics {\bf 102} (1975),  223-290.

 \smallskip

 \noindent 
\noindent
[HL$_2$] F. R. Harvey and H. B. Lawson, Jr, {\sl On boundaries of complex
analytic varieties, II}, Annals of Mathematics {\bf 106} (1977), 
213-238.

 \smallskip

\noindent
[HL$_3$] F. R. Harvey and H. B. Lawson, Jr, {\sl An intrinsic
characterization of K\"ahler manifolds}, Inventiones Math.,  {\bf 74} 
(1983), 169-198.
 \smallskip

 \noindent
[HL$_4$] F. R. Harvey and H. B. Lawson, Jr, {\sl Boundaries of 
varieties in projective manifolds}, J. Geom. Analysis,  {\bf 14}
no. 4 (2005), 673-695. ArXiv:math.CV/0512490.
 \smallskip

 \noindent 
\noindent
[HL$_5$] F. R. Harvey and H. B. Lawson, Jr, {\sl Projective hulls and
the projective Gelfand transformation}, Asian J. Math. {\bf 10}, no. 2 (2006), 279-318. ArXiv:math.CV/0510286.

 \smallskip

\noindent
[HL$_6$] F. R. Harvey and H. B. Lawson, Jr, {\sl Remarks on the Alexander-Wermer Theorem}, 
ArXiv: math.CV/0610611.

 \smallskip

\noindent
[HL$_7$] F. R. Harvey and H. B. Lawson, Jr, {\sl Projective linking and boundaries of holomorphic chains in projective manifolds, Part II},  pp. 365-380 in ``Inspired  by S. S. Chern'', P. Griffiths, Ed., 
Nankai Tracts in Math. {\bf 11}, World Scientific Publishing, Singapore, 2006.

 \smallskip

\noindent
[HL$_8$] F. R. Harvey and H. B. Lawson, Jr, {\sl Relative holomorphic cycles and duality}, Preprint, Stony Brook, 2006.

 \smallskip

\noindent
[HLZ] F. R. Harvey, H. B. Lawson, Jr. and J. Zweck, {\sl A
deRham-Federer theory of differential characters and character duality},
Amer. J. of Math.  {\bf 125} (2003), 791-847.

 \smallskip

\noindent
[HLW] F. R. Harvey, H. B. Lawson, Jr. and J. Wermer, {\sl 
The projective hull of certain curves in $\bbc^2$}, ArXiv:math.CV/0611482.

 \smallskip

\noindent
[HP] F. R. Harvey, J. Polking, {\sl Extending analytic objects},
Comm.  Pure Appl. Math. {\bf 28} (1975), 701-727.

 \smallskip

   \noindent
[S]   H. H. Schaefer,  Topological Vector Spaces,    Springer Verlag,
New York,  1999.

\smallskip

 \noindent
[W$_1$]   J. Wermer  {\sl    The hull of a curve in $\bbc^n$},    
Ann. of Math., {\bf  68}  (1958), 550-561.

\smallskip

 \noindent
[W$_2$]   J. Wermer    {\sl    The argument principle and
boundaries of analytic varieties},     Operator Theory: Advances and 
Applications, {\bf  127}  (2001), 639-659.

\smallskip

\end